\documentclass[a4paper,12pt]{article}
\usepackage[utf8]{inputenc}
\usepackage[T1]{fontenc}
\usepackage[english]{babel}
\usepackage{amsmath,amsfonts,dsfont}
\usepackage{mathabx}
\usepackage{mathrsfs}
\usepackage{enumerate}
\usepackage{fullpage}
\usepackage{lscape}
\usepackage{stmaryrd}
\SetSymbolFont{stmry}{bold}{U}{stmry}{m}{n}
\usepackage{graphicx}
\usepackage[clockwise]{rotating}
\usepackage{silence}
\usepackage{url}

\newcommand{\M}{\mathcal{M}}
\newcommand{\J}{\textbf{J}}

\newcommand{\Mh}{\mathcal{M}_h}

\renewcommand{\phi}{\varphi}
\renewcommand{\epsilon}{\varepsilon}

\newcommand{\x}{\boldsymbol{x}}
\newcommand{\xn}{\boldsymbol{\vec{x}\cdot\n}}

\newcommand{\s}{\boldsymbol{s}}

\newcommand{\W}{{\textbf W}}

\newcommand{\F}{\vec{\boldsymbol{\mathcal{F}}}}
\newcommand{\FStar}{\textbf{F}^*}

\newcommand{\V}{{\bf V}}
\newcommand{\X}{{\textbf{X}}}

\newcommand{\demi}{\frac{1}{2}}

\newcommand{\n}{\vec{\bf n}}

\renewcommand{\P}{\text{P}}

\renewcommand{\S}{\boldsymbol{\mathcal{S}}}

\renewcommand{\O}{\mathcal{O}}

\renewcommand{\L}{\text{L}}

\newcommand{\aire}[1]{\lvert #1 \rvert}

\newcommand{\G}{\mathscr{G}}

\newcommand{\K}{\mathcal{K}}

\newcommand{\saut}[1]{\llbracket #1 \rrbracket}

\newcommand{\abs}[1]{\left| #1 \right|}

\newcommand{\ds}{\displaystyle}

\newcommand{\Vect}[1]{\overrightarrow{#1}}
\newcommand{\Div}[1]{\Vect{\nabla}\cdot #1}

\newcommand{\N}{\mathbb{N}}

\newcommand{\QQ}{\mathbb{Q}}
\newcommand{\One}{\mathds{1}}
\newcommand{\R}{\mathbb{R}}

\newcommand{\T}{\text{T}}

\newcommand{\Q}{\text{Q}}

\renewcommand{\emph}{\textit}

\newcommand{\beq}[1]{\begin{equation}\label{#1}}
\newcommand{\eeq}{\end{equation}}
\newcommand{\bitem}{\begin{itemize}}
\newcommand{\eitem}{\end{itemize}}
\newcommand{\benum}{\begin{enumerate}}
\newcommand{\eenum}{\end{enumerate}}

\usepackage[L]{thmbox}
\newtheorem[M]{remark}{Remark}
\newtheorem[M]{example}{Example}
\newtheorem[M]{homework}[remark]{Homework}
\newtheorem[M]{property}[remark]{Property}
\newtheorem[L]{principle}[remark]{Principle}
\newtheorem[L]{theorem}[remark]{Theorem}
\newtheorem[L]{lemma}[remark]{Lemma}
\newtheorem[L]{corrolary}[remark]{Corrolary}
\newtheorem[M]{definition}[remark]{Definition} 

\usepackage{tikz}
\usetikzlibrary{shapes}
\tikzstyle{mybox} = [draw=black, fill=white, thick,
    rectangle, rounded corners, inner sep=10pt, inner ysep=20pt]
\tikzstyle{fancytitle} =[fill=black, text=white]
\usepackage{xargs}
\newcommandx{\fancybox}[3][1=Title of the box, 2=0.99\textwidth]{%
\begin{center}
\begin{tikzpicture}
\node [mybox] (box){%
    \begin{minipage}{#2}
    #3
    \end{minipage}
};
\node[fancytitle, right=10pt] at (box.north west) {#1};
\end{tikzpicture}%
\end{center}
}
\newcommand*\circled[1]{\tikz[baseline=(char.base)]{
            \node[shape=circle,draw,inner sep=2pt] (char) {#1};}}

\title{Notes on the Discontinuous Galerkin methods for the numerical simulation of hyperbolic equations}
\author{Adam Larat}
\date{\today}

\begin{document}

\maketitle

\section{General Context}

\subsection{Bibliography}

The roots of Discontinuous Galerkin (DG) methods is usually attributed to Reed and Hills in a paper published in 1973 on the numerical approximation of the neutron transport equation \cite{ReedHill73}.
In fact, the adventure really started with a rather thoroughfull series of five papers by Cockburn and Shu in the late 80's \cite{CockburnShu1,CockburnShu2,CockburnShu3,CockburnShu4,CockburnShu5}. 
Then, the fame of the method, which could be seen as a compromise between Finite Elements (the center of the method being a weak formulation) 
and Finite Volumes (the basis functions are defined cell-wise, the cells being the elements of the primal mesh) increased and slowly investigated 
successfully all the domains of Partial Differential Equations numerical integration. In particular, one can cite the ground papers for the common treatment of 
convection-diffusion equations \cite{BassiRebay1,BassiRebay2} or the treatment of pure elliptic equations \cite{Arnold01,OrtnerSuli07}. 
For more information on the history of Discontinuous Galerkin method, please refer to section 1.1 of \cite{HesthavenWarburton}.

Today, DG methods are widely used in all kind of manners and have applications in almost all fields of applied mathematics. (TODO: cite applications and 
structured/unstructured meshes, steady/unsteady, etc...). The methods is now mature enough to deserve entire text books, among which I cite
a reference book on Nodal DG Methods by Henthaven and Warburton \cite{HesthavenWarburton} with the ground basis of DG integration, numerical analysis of its linear behavior 
and generalization to multiple dimensions. 

Lately, since 2010, thanks to a ground work of Zhang and Shu \cite{ZhangShu10a,ZhangShu10b,ZhangThesis,ZhangShu11,ZhangXiaShu12}, Discontinuous Galerkin methods are eventually able to combine high order accuracy 
and certain preservation of convex constraints, such as the positivity of a given quantity, for example. These new steps forward are very promising since it brings us 
very close to the "\textit{Ultimate Conservative Scheme}", \cite{VanLeer79,Abgrall01}. 

\subsection{Hyperbolic conservation laws}

This section is just a quick introduction to hyperbolic conservation laws in order to settle the notations. 

Let $d\in\N^*$ be the number of spatial dimensions and $\Omega\subset\R^d$ be the domain of study. $\partial\Omega$ denotes the frontiers of the domain. 
To greatly simplify the following notes, we will ignore boundary conditions. 

Let $m\in\N^*$ be the number of conserved variables and $\W$ the vector of these variables. Then, $\W$ is ruled by a \textit{Conservation Law} if there 
exists a \textit{Flux Function}
\begin{equation}\label{eq:FluxFunction}
\left\{
  \F:
  \begin{array}{ccc}
    \S\subset\R^m & \longrightarrow & \left(\R^m\right)^d\\
    \W & \longmapsto & \F(\W)
  \end{array}
\right.
\end{equation}
such that $\W$ verifies the PDE: 
\beq{eq:CL}
  \frac{\partial \W}{\partial t} + \Div{\F\left(\W\right)} = 0.
\eeq
Here, the set $\S$ denotes the possible physical constraints on the conserved variables (positivity of the density, 
realizability of the transported moments, etc\dots). 

This conservation law is said to be \textit{hyperbolic}, if the Jacobian 
\beq{eq:Jacobian} \J\left(\W\right) = \dfrac{\partial\F}{\partial\W}\eeq 
is diagonalizable with real eigenvalues. In this case, we denote $\alpha$ the maximal absolute eigenvalue of $\J$:
\beq{eq:Alpha} 
  \alpha\left(\W\right) 
  = 
  \max \left\{\abs{\lambda},\; \exists\V\in\R^m \text{ such that } \J\left(\W\right).\V = \lambda \V\right\}.
\eeq

\section{Runge-Kutta discontinuous Galerkin methods}

\subsection{Degrees of freedom}\label{ssec:DoFs}

Let $\Omega\subset\R^d$ be the \textit{domain of study} and 
\beq{eq:Tesselation} 
  \bigcup_{k=1}^{N} \T_k = \Omega 
  \quad \text{ such that } \quad 
  \forall i,j \in \saut{1,N},\quad 
  \text{dim}\left(\T_i\cap\T_j\right)<d,
\eeq
a \textit{tesselation} of the domain. 

In each element $\T_k$, we define a funcitonal basis $\ds \left(\psi_l^k\right)_{l=1,\dots,\P_k}$ and associate $\T_k$ with the finite dimensional 
vector space:
\beq{eq:functionalSpace}
  V_h^k = \left\{
    \phi:\Omega \longrightarrow \R,
    \quad
    \phi = \left(\sum_{l=1}^{\P_k} a_l \psi_l^k\right)\chi_k
    \right\},
\eeq
where $\chi_k$ denote the characteristic function of $\T_k$, taking value $1$ in $\T_k$ and $0$ elsewhere. 

Therefore, a numerical solution on $\M_h = \left(\T_k\right)_{k=1,\dots,N}$ has exactly 
$\ds \P = \sum_{k=1}^{N} \P_k$
\textit{Degrees of Freedom (DoF)}.

\subsection{Spatial weak formulation}

The continuous hyperbolic conservation law \eqref{eq:CL} is approximated on the high order functional vectorial space $V_h$ defined 
in \eqref{eq:functionalSpace}. It reads: 
\fancybox[Weak Formulation]{
  Find $\ds \W_h(t,\x) = \sum_{k=1}^{N} \chi_k(\x) \sum_{j=1}^{\P_k} \W_k^j(t) \psi_k^j(\x)$ a piecewise 
  high order polynomial solution, such that:
  \beq{eq:WeakFormulation}
    \forall k' \in\saut{1,N}, \forall i \in \saut{1,\P_{k'}}, \quad 
    \int_{\Omega} \left(\partial_t \W_h(t,\x) + \Div{\F\left(\W_h(t,\x)\right)}\right) \psi_{k'}^i(\x) d\x = 0.
  \eeq
}

Since the supports $\T_k$ and $\T_{k'}$ are disjointed, $k'=k$ always and the problem comes to 
$$
\forall k, i, \quad 
\int_{\T_k} \sum_j d_t\W_k^j(t) \psi_k^j(\x)\psi_k^i(\x) d\x + \int_{\T_k} \Div{\F\left(\W_h(t,\x)\right)}\psi_k^i(\x)d\x = 0.
$$
By integration by parts, one easily gets the following three fold formulation: 
\beq{eq:DGNum}
  \underbrace{\sum_j\left(\int_{\T_k} \psi_k^j\psi_k^i d\x\right)\, d_t\W_k^j(t)}_{\circled{1}}
  - 
  \underbrace{\int_{\T_k} \F\left(\W_h(t,\x)\right)\cdot\nabla\psi_k^i(\x)d\x}_{\circled{2}}
  +
  \underbrace{\int_{\partial\T_k} \psi_k^i(\s)\,\F\left(\W_h(t,\s)\right)\cdot\n\, d\s}_{\circled{3}}
  = 0
\eeq

\begin{remark}
  One usually considers (but we don't have to)
  \beq{eq:FluxProjection}
    \F\left(\W_h(t,\x)\right) = \sum_j\F\left(\W_k^j(t)\right)\psi_k^j(\x).
  \eeq
  If not, \circled{2} is integrated by means of an adequate volumic quadrature. 
\end{remark}

\begin{remark}
  Along the boundaries $\partial\T_k$ of $\T_k$, $\F\left(\W_h\right)$ does not have a mathematical sens 
  other than the flux occuring locally between the two states on both sides. Then the flux appearing in 
  \circled{3} is replaced by the numerical flux: 
  \beq{eq:NumFlux}
    \F\left(\W_h(t,\s)\right)\cdot\n = \FStar\left(\W_{\text{ext}}(t,\s),\W_{\text{int}}(t,\s);\n\right).
  \eeq
\end{remark}

In the case the flux function is also projected on the basis as in \eqref{eq:FluxProjection}, the global formulation 
simply comes down to:
\beq{eq:FinalDG}
  \M_{ij}d_t\W_j 
  - 
  \vec{\K}_{ij}^t\cdot\F\left(\W_j\right) 
  +
  \sum_{\epsilon\in\partial\T_k} 
  \underbrace{\int_{\epsilon} \psi_k^i(\s)  \FStar\left(\W_{\text{ext}}(t,\s),\W_{\text{int}}(t,\s);\n\right) d\s}_{\circled{4}}
  = 
  0
\eeq

\begin{remark}
  \bitem
    \item In general, last integral \circled{4} is estimated by mean of an adequate quadrature of dimension $d-1$. 
       This implies a numerous amount of numerical flux computations. 
    \item Luckily, since we transport and update all the derivatives of the solution within the mesh, we don't 
      need to be high-order accurate on all these numerical fluxes. 
      A first order flux on all the degrees of freedom will provide a high order 
      solution anyway. Therefore, the linear Rusanov flux \eqref{eq:Rusanov} is very well suited for this purpose, since 
      its evaluation is quite costless and it has good monotonicity properties (see section \ref{sssec:PositiveFluxes}). 
  \eitem
\end{remark}

\subsection{Basis functions}\label{ssec:BasisFunctions}

DG methods are split into two main families: \textit{modal} and \textit{nodal} DG. These prefixes refer 
to the choice on the basis functions $\psi_l^k$. In general, the modal basis function will be chosen so 
that the corresponding mass matrix
\beq{eq:MassMatrix}
\M_{ij}^k = \int_{\T_k} \psi_i^k \psi_j^k 
\eeq 
is everywhere diagonal, what greatly simplifies the updates. On the other hand, any access to the value of the 
solution at a given point (for example for the evaluation of a numerical flux somewhere) will require a combination of all 
the degrees of freedom in the element: 
\beq{eq:PointEval}
  \W_h(x_k) = \sum_{l=1}^{\P_k} a_l \psi_l^k(x_k).
\eeq

On the contrary, nodal basis function are attached to certain points of the considered element. In general, the points 
are chosen in a way it simplifies the volumic and edge terms evaluations. In particular, on the edges of the elements degrees of 
freedom should be attached to some known quadrature points. Then, the numerical flux integral evaluation has all the necessary data
for the numerous numerical fluxes evaluation directly at hand. 

For dimensions stricly higher than 1 and polynomial order at least quadratic, there is a conjecture that there is no orthogonal 
nodal basis\footnote{Personal conversation with V. Perrier}. This might even be proven, alas I don't know yet about the demonstration. 

\subsubsection{Modal basis}

\bitem
\item \textbf{Legendre polynomials:} this is the orthogonal interlocked polynomial basis for measure one. It is simply obtained by 
a Gram-Schmidt orthogonalization process starting from any polynomial basis of increasing order. 
\beq{eq:Legendre1D}
\text{On } [-1,1]: \quad 1, \; x,\; x^2-1/3,\, \dots
\eeq
\item Other orthogonal basis: Jacobi, Tchebychev, trigonometric, etc.
\eitem 

\subsubsection{Nodal basis}

Let us restrict to 1D for the moment. Givent a set $\left(x_i\right)_{i=1,\dots,\P_k}$ of points in $\T_k$, the associated \textit{Lagrangian} basis
is the only basis of polynomials of order $\P_k -1$ such that $\L_j(x_i) = \delta_{ij}$. Then it is easy and classic to write 
\beq{eq:Lagrange1D}
  \L_i(x) = \dfrac{\ds\prod_{j\neq i} (x-x_j)}{\ds\prod_{j\neq i} (x_i -x_j)}.
\eeq

In general, the set of nodes is: 
\bitem
  \item a regular distribution within the cell: 
    \beq{eq:RegDistribution} x_j = x_{i-1/2} + j \frac{\Delta x}{\P_k-1}, \quad j = 0,\dots,\P_k-1,\eeq
  \item a set of smart quadrature points, like Gauss-Lobatto quadrature points since they contain the bounds of the interval. 
\eitem

In more dimensions, the generalization of the former discussion is more complex. Only the regular distribution can be algorithmically extended. 
In particular, the generalization of the Gauss-Lobatto Lagrange polynomials to 2 and 3 D is not obvious, see discution in \cite{HesthavenWarburton}, chapter 6. 

\subsection{Time integration}

Now, thanks to the spatial semi-integration, the numerical scheme \eqref{eq:FinalDG} comes down to an Ordinary Differential Equation.
It can be integrated by any numerical procedure for ODEs. Here we focus on Runge-Kutta methods, since they offer a hierarchy 
of increasing order, which is suitable for the time integration of high order spatial weak formulations. 

\subsubsection{Runge-Kutta methods}

Let \beq{eq:ODE} y'(t) = f(t,y(t)) \eeq be an arbitrary ordinary differential equation. 
An $s$-stages explicit Runge-Kutta method applied to this equation can be written in the form:
\beq{eq:RKExp}
  \left\{
    \begin{array}{ccl}
      y_{n+1} &=& \ds y_n + \Delta t \sum_{i=1}^s b_i k_i,\\
      
      k_i     &=& \ds f\left(t_n+c_i\Delta t, y_n +\Delta t \sum_{j<i}a_{ij}k_j\right).
    \end{array}
  \right.
\eeq
This writing can be summed up in the so-called \textit{Butcher Tableau}: 
\beq{eq:Butcher}
  \begin{array}{c|cccc}
    0 & & 0 && \\
    c_2& a_{21}&&& \\
    \vdots&\vdots&\ddots&&\\
    c_s&a_{s1}&\cdots&a_{s,s-1}&\\
    \hline
       &b_1&\cdots&\cdots&b_s
  \end{array}
\eeq

These methods are consistent with equation \eqref{eq:ODE} when 
\beq{eq:consistency} c_i = \sum_j a_{ij}. \eeq
Then, under the additional constraint that $\sum b_j = 1$, the explicit numerical procedure 
is stable if the time step is smaller than a bound which depends on $f$ Lipschitz constant. 

In certain cases, especially when dealing with stiff multiscale problems, 
this stability constraint may be too harsh. One can overcome it by going implicit. The intermediate updates now become 
\beq{eq:RKImp}
  k_i = \ds f\left(t_n+c_i\Delta t, y_n +\Delta t \sum_{j}a_{ij}k_j\right),
\eeq
and the Butcher tableau now looks like
\beq{eq:ButcherImp}
  \begin{array}{c|ccc}
    c_1&&& \\
    \vdots&&A&\\
    c_s&&&\\
    \hline
       &b_1&\cdots&b_s
  \end{array}
\eeq
where $A$ is a full matrix. 
The price to pay for the increased range of stability is a interdependance between all the intermediate steps, which 
is not the case in the explicit version, thanks to the triangular shape of $A$. The solution at time $t+\Delta t$ is 
usually the solution of a big non-linear system. 

The accuracy study of Runge-Kutta methods is more complex. All I want to say here is that it is known that starting from 
$s\geq 5$, there is no more RK method of order $s$. 

\subsubsection{Strong Stability Preserving integrators}

In a fundamental paper published in 1988 \cite{ShuSSP88}, Chi-Wang Shu selects among all the RK methods a family  
of integrators he first calls "\textit{Total Variation Diminishing time discretizations}", that will later be renamed as
\textit{Strong Stability Preserving} time discretizations in a review paper \cite{Gottlieb01}. 

These methods can be seen as those which can be written as a convex combination of Explicit Euler (EE) time integrations:
\beq{eq:SSP}
  y_{n+1} = \sum_{i=0}^{s-1} a_i y_{(i)} + \Delta t b_i f\left(t+i\Delta t,y_{(i)}\right).
\eeq
Since the stability constraint of the EE scheme is $1.0$, these methods are stable under the fact that the time step 
$\Delta t$ is smaller than the Lipschitz constant of $f$ time 
\beq{eq:CFL_SSP}
  \text{CFL}_{\text{max}} = \max_i \frac{a_i}{b_i}.
\eeq

Next, Shu and Osher \cite{ShuOsher88,ShuOsherII89} proved that up to fourth order, there exist an optimal SSP-RK integrator, meaning 
being both
\bitem
  \item of the order of the number of stages, 
  \item with the largest possible stability constraint: $$\text{CFL}_{\text{max}} = 1.0.$$
\eitem

\begin{example}[Heun's Method]
This is the optimal second order SSP-RK integration : 
\beq{eq:Heun}
  \left\{
    \begin{array}{ccc}
      y_{n+\demi} &=& y_n + \Delta t f\left(t,y_n\right),\\
      y_{n+1}     &=& \demi\left[y_n+\left(y_{n+\demi}+\Delta t f\left(t+\Delta t, y_{n+\demi}\right)\right)\right].
    \end{array}
  \right.
\eeq
\end{example}

\subsection{Convex state preserving DG methods}

\subsubsection{Convex state preserving numerical flux}\label{sssec:PositiveFluxes}

Let $\S\subset\R^m$, convex, be the set of physical states. 
\begin{example}[Convex Constraints]
\bitem 
  \item \textbf{Euler Equations:} Density and pressure have to stay positive. If the Equation Of State (EOS) follows the 
    Bethe-Weil conditions, these conditions imply the convexity of the physical states\footnote{Not shure about that last remark. 
    Needs to be checked. If the EOS is the perfect gas law, it is ok.}.
  \item \textbf{Moments of a Repartition Function:} if the transported conserved quantities are the moments of a positive 
    repartition function, then these moments need to stay moments of a positive distribution. In many cases, this implies 
    convex constraints on the set of moments, \cite{DetteStudden}.
\eitem 
\end{example}

\fancybox[CSP Numerical Flux] {
  A numerical flux $\ds\FStar\left(\W^+,\W^-;\n\right)$ is said to be \textbf{Convex State Preserving (CSP)}, when, 
  for any states
    $$ \W_{i-1}^n, \W_i^n, \W_{i+1}^n \in \S,$$
  the Explicit Euler update belongs to $\S$,
    \beq{eq:EEUpdate} 
      \W_i^{n+1} 
      = 
      \W_i^n - \nu \left[\FStar\left(\W_{i+1}^n,\W_i^n;\n\right)-\FStar\left(\W_i^n,\W_{i-1}^n;\n\right)\right] 
      \in\S,
    \eeq
  under a CFL constraint:
    \beq{eq:CFLConstraint}
      \nu\leq C.
    \eeq
}

\begin{example}
  \bitem
    \item \textbf{Rusanov Flux:} 
      \beq{eq:Rusanov} 
        \FStar\left(\W^+,\W^-;\n\right) = 
        \frac{\F(\W^+)+\F(\W^-)}{2}\cdot\n - \alpha\left(\W^+-\W^-\right).
      \eeq
      $\alpha$ being defined by \eqref{eq:Alpha}. 
    \item \textbf{Godunov Flux:} since the Godunov flux is a convex combination of the physical states coming from the exact 
      resolution of the Riemann problem at the interface, the update is physical and \eqref{eq:EEUpdate} comes true. 
    \item \textbf{HLLx Solvers:} the same reasoning applies to HLL solvers, since the updated states are convex combinations 
      of physical states, even though the resolution of the Riemann problem is approximated. 
   \eitem
\end{example}

\begin{corrolary}[SSP Methods]\label{Corr:SSP}
  This naturally extends to SSP integrators, since they are convex combination of Explicit Euler updates. Only 
  the CFL condition must be multiplied by the additional constraint \eqref{eq:CFL_SSP}. Hence the 
  particular role of optimal SSP-RK integrators.
\end{corrolary}

\subsubsection{Application to RK-DG methods}

Recall the general update of all the degrees of freedom of an element $\T_k$: 
\begin{equation*}\tag{\ref{eq:FinalDG}}
  \M_{ij}d_t\W_j^k 
  - 
  \int_{\T_k} \F\left(\W(t,\x)\right)\cdot\nabla\psi_i^k \; d\x 
  +
  \int_{\partial\T_k} \psi_k^i(\s)  \FStar\left(\W_{\text{ext}}(t,\s),\W_{\text{int}}(t,\s);\n\right)\;d\s
  = 
  0.
\end{equation*}
Once summed up on all the DoFs, we obtain the equation ruling the \textit{mean value} $\overline{\W_k}$:
\beq{eq:MeanValue}
  \aire{\T_k}\frac{d\overline{\W_k}}{dt} + \int_{\partial\T_k} \FStar\cdot\n\;d\s  = 0.
\eeq
If the ODE \eqref{eq:FinalDG} is integrated by Explicit Euler, the update simply gives:
\beq{eq:MeanValueDiscrete}
  \overline{\W_k}^{n+1} =  \overline{\W_k}^{n} - \frac{\Delta t}{\aire{\T_k}}\int_{\partial\T_k} \FStar\cdot\n\;d\s  = 0.
\eeq
Note that by Corrolary \ref{Corr:SSP}, this discussion naturally extends to SSP integrators. 

Now, let's switch to 1D and assume that the polynomial order is $\P_k-1$, so that there exists a number $Q$ of 
Gauss-Lobatto quadrature points, $(2Q-3)\geq (\P_k-1)$, such that the following equality is exact:
\beq{eq:GLExact}
  \overline{\W_k} = \sum_{q=1}^Q \omega_q \W_k(x_q).
\eeq
Then, discrete update \eqref{eq:MeanValueDiscrete} can be rewritten into
\begin{eqnarray}
\overline{\W_k}^{n+1} 
&=& \sum_{q=1}^Q \omega_q \W_k^n(x_q) - \frac{\Delta t}{\Delta x_k} \bigg\{ \FStar\left(\W_{k+1}^-,\W_k^+\right)-\FStar\left(\W_k^-,\W_{k-1}^+\right)\bigg\}, \label{eq:ConvexUpdate}\\
&=& \sum_{q=1}^Q \omega_q \left[\W_k^n(x_q) - \frac{\Delta t}{\omega_q\Delta x_k}
  \bigg\{ \FStar\Big(\W_k^n(x_{q+1}),\W_k^n(x_{q})\Big)-\FStar\Big(\W_k^n(x_q),\W_k^n(x_{q-1})\Big)\bigg\}\right],
  \nonumber
\end{eqnarray}
with the convention that $\W_k^n(x_{Q+1}) = \W_{k+1}^-$ and $\W_k^n(x_0) = \W_{k-1}^+$. 
In this form, we see that the DG update in the mean writes as a convex combination of abstract Euler Explicit updates at the 
Gauss-Lobatto quadrature points. So that
\bitem 
  \item if the numerical flux is CSP \eqref{eq:EEUpdate},
  \item  if the solution at time $t_n$ belongs to $\S$ at all the Gauss-Lobatto quadrature 
    points, 
    \beq{eq:WatGL}
      \forall q\in\saut{0,Q+1}, \quad \W_k^n(x_q)\in\S,
    \eeq
\eitem
the updated mean value $\overline{\W_k}^{n+1}$ will be in $\S$, under a CFL constraint
\beq{eq:CFLConvexConstraint}
  \nu \leq \omega_1 C, 
\eeq
since it is known that the smallest weights of the Gauss-Lobatto quadrature are always at the borders of the interval: 
$\omega_1 = \omega_Q = \min \omega_q$.

Starting from that, Zhang and Shu \cite{ZhangShu10a} provided a limitation procedure which was proven to conserve accuracy. 
Given that the updated mean value $\overline{\W_k}^{n+1}$ belongs to $\S$ and that $\S$ is a convex set, we now look 
at the values of the update solution at the Gauss-Lobatto quadrature points. If any of these values $\W_k^{n+1}(x_q)$ 
is outside $\S$, there exists a unique value $\theta_q\in]0,1[$, 
such that $\theta_q \W_k^{n+1}(x_q) + (1-\theta_q) \overline{\W_k}^{n+1}$ is back on the frontier of $\S$.
By setting
\beq{eq:theta}
  \theta = \min_{q} \theta_q, 
\eeq
and 
\beq{eq:LimitedW}
  \widetilde{\W}_k^{n+1} = \theta \W_k^{n+1} + (1-\theta)\overline{\W_k}^{n+1},
\eeq
one obtains a new solution which is as accurate as $\W_k^{n+1}$ but limited in a way that \eqref{eq:ConvexUpdate}
will propagate the convex constraint preservation further. 

\begin{remark}
\bitem
  \item In a paper published on Arxiv in 2012, \cite{JohnsonRossmanith13}, Johnson and Rossmanith proved that such 
    a limitation procedure on Gauss-Lobatto quadrature points is optimal in 1D. 
  \item We have restrict the discussion to 1D because it is much simpler to explain in this context. However, 
    the procedure generalizes to any dimension \cite{ZhangShu10b,ZhangXiaShu12,JohnsonRossmanith13} but in a 
    rather cumbersome manner. The interested reader is encouraged to read these articles and citation within. 
\eitem 
\end{remark}

\subsection{Limits of such procedure}

Looking for a Strong Stability Preserving strategy or not, generalization of Runge-Kutta integration to higher orders is 
rather complex, since an increasing number of intermediate stages is needed to obtain an additional order in time. 
Moreover, all these intermediate stages generally need to be stored 
for the final $t^{n+1}$ update. This may be limiting in term of memory usage. 

On the top of that, the global time step of the method is chosen a priori and kept during the whole multi-step Runge-Kutta process. Even if a 
security coefficient is applied globally on the time step, the local physics may evolve rapidly within the time step and the stability 
constraint on the time progress may be a posteriori violated. This risk increasing with the number of substeps. 

These reasons explain why people start to turn to space-time formulation. A fully space-time formulation of DG methods is conceivable but 
finally looks like a huge implicit formulation: all the degrees of freedom of the mesh are coupled, especially with non-linear equations.

On the other hand, arbitrary high order one-step explicit space-time have appeared under the name of ADER (for \textbf{A}rbitrary high-order schemes using \textbf{DER}ivatives).
Even though these formulations are intrinsically space-time, the volume and flux terms \circled{2} and \circled{3}, in \eqref{eq:DGNum}, are in fact somehow 
extrapolated from the available information at time $t^n$ and the method can be considered as one-step and explicit.
This is the topic of the next section. 

\section{ADER-DG integration}

ADER-DG methods always start by a space-time integration of a weak formulation in space of equation \eqref{eq:CL}. 
\begin{eqnarray}
  & &
  \int_{t^n}^{t^{n+1}} 
  \left[ 
    \int_{\T_k}\left(\frac{\partial \W}{\partial t} + \Div{\F\left(\W\right)}\right)
               \cdot
               \psi_i^k(\X)\, d\x
  \right]
  \, dt = 0,\label{eq:ADERDG}\\
  &\Leftrightarrow&\nonumber\\
  & &\M_{ij} \left(\W_j^{k,n+1}-\W_j^{k,n}\right) 
  - 
  \underbrace{\int_{t^n}^{t^{n+1}}\int_{\T_k} \F\left(\W\right)\cdot\nabla\psi_i^k \, d\x dt}_{\circled{5}} 
  + 
  \underbrace{\int_{t^n}^{t^{n+1}}\int_{\partial\T_k}\psi_i^k\FStar(t,\n) \, d\s \,dt}_{\circled{6}}
  = 
  0.\nonumber
\end{eqnarray}
In the two next sections, the goal is to set up a method to obtain the necessary data needed to integrate 
the two integrals \circled{5} and \circled{6} in last equation. 
In a first procedure, the time evolution is gotten from an arbitrary high order 
Taylor expansion in time of the conserved variables. Thanks to the \textit{Cauchy-Kowaleski} procedure 
(also known as \textit{Lax-Wendroff} procedure), the time derivatives of $\W$ are functions of the spatial derivatives of the 
fluxes, which are known at time $t^n$. A \textit{Generalized Riemann Problem} (GRP) is solved at each interface 
and the problem can be moved to next time step. 

In a second procedure, a predictor step is ran in the form of a \textit{local space-time DG scheme}. This 
allows to get a local information on the evolution of the data. The predicted solution is then used to 
compute the space-time flux and volume terms of equation \eqref{eq:ADERDG}. 

\subsection{Cauchy-Kowaleski ADER procedure}

\subsubsection{Generalized Riemann Problem}

Most of this paragraph takes its roots in chapter 19 and 20 of Toro's book \cite{toro1999}.

In a general context, the numerical flux needed to compute last integral \circled{6} of equation \eqref{eq:ADERDG} 
comes from the (possibly approximate) resolution of a Riemann problem:
\beq{eq:RP}
  \left\{
    \begin{array}{l}
      \ds \partial_t \W + \partial_{\n}\F\left(\W\right) = 0, \quad t>0,\; \xn\in\R,\\
      \\
      \ds \W(t=0,\xn<0) = \W^-, \\
      \ds \W(t=0,\xn>0) = \W^+.
    \end{array}
  \right.
\eeq
By autosimilarity, this solution depends only on one variable $\xi=\xn/t$ and the numerical flux usually writes
\beq{eq:NumFluxRP}
  \FStar(\W^+,\W^-;\n) = \F\left(\W^*(\xi=0)\right)\cdot\n,
\eeq
where $\W^*(\xi=0)$ is the values taken by the (possibly approximate) solution of the Riemann Problem along 
the ordinate axis. 

At higher order in space, input values $\W^+$ and $\W^-$ are fed with the limit values at the interface 
of the polynomial in each cell: $\W_{k+1}^n(x_{i+\demi}^+)$ and $\W_{k}^n(x_{i+\demi}^-)$.
Therefore a first order error in time comes from the fact that the Riemann solvers see a constant extrapolation 
of $\W^+$ and $\W^-$ in the neighboring cells. 

This stimulated people to look at \textit{Generalized Riemann Problems} (GRP)
where the numerical flux would now be:
\beq{eq:NumFluxGRP}
  \FStar\left(t;\W^+(x>0),\W^-(x<0);\n\right) = \F\left(\W^*(t,x=0)\right)\cdot\n, \quad t\in]0,\Delta t[.
\eeq
The accuracy in time would then be achieved by a Taylor expansion of $\W$ in time along the ordinate axis
\beq{eq:TaylorExpansion}
  \W^*(t,x=0) = \W^*(0^+,x=0) + \sum_{k=1}^N\frac{t^k}{k!}\frac{\partial^k\W^*}{\partial t^k}(0^+,x=0) + \O\left(t^{N+1}\right),
\eeq
and to look at an accurate enough way to solve all the unknowns (meaning the successive time derivatives at $t=0^+$). 
It came out that the problem can be split into $N$ sub-Riemann problems, one being the full space consistent Riemann problem
\eqref{eq:RP}, called the \textit{space high-order Riemann problem}, the $(N-1)$ others being rather simple.

In fact, since $\W$ is the solution of \eqref{eq:CL}, its $k^{\text{th}}$-order time derivatives can be written as a function 
of its spatial derivatives up to order $k$: 
\beq{eq:LWExpansion}
  \partial_t^{(k)}\W(t,x) = \G\left(\W(t,x),\partial_x^{(1)}\W(t,x),\dots,\partial_x^{(k)}\W(t,x)\right).
\eeq
Now, one can derivate problem \eqref{eq:RP} $k$ times in space to get a Riemann problem on $\partial_x^{(k)}\W(t,x)$.
In fact, since only the solution at $(t=0^+,x=0)$ is needed, the associated Riemann problem can be simplified into 
\beq{eq:RPk}
  \left\{
    \begin{array}{l}
      \partial_t\left(\partial_x^{(k)}\W(t,x)\right) 
      + 
      \J\Big(\W^*(t=0^+,x=0)\Big)\partial_x\left(\partial_x^{(k)}\W(t,x)\right) 
      = 0,\\
      \\
      \partial_x^{(k)}\W(t=0,x) = \left\{
        \begin{array}{l}
          \partial_x^{(k)}\W(t=0,x=0^-) \; \text{ if } \; x<0, \\
          \partial_x^{(k)}\W(t=0,x=0^+) \; \text{ if } \; x>0,
        \end{array}
      \right.
    \end{array}
  \right.
\eeq
where $\W^*(t=0^+,x=0)$ is the solution of the space high-order Riemann problem and $\J$ is the Jacobian of the flux 
at this state.

\subsubsection{Integration within DG formulation}

Once this is done, the Taylor expansion \eqref{eq:TaylorExpansion} can be filled up to desired order $N$ 
thanks to the Lax-Wendroff expansions \eqref{eq:LWExpansion}, and the numerical flux is gained as the exact integration 
 on $t\in[0,\Delta t]$ of 
 $$ \F\Big(\W^*(t,x=0)\Big).$$ 

Next, the volume integral \circled{5} needs to be evaluated. This is generally done by a numerical quadrature 
at a sufficient order of accuracy:
\beq{eq:XTQuadrature}
  \int_{t^n}^{t^{n+1}}\int_{\T_k} \F\left(\W(t,\x)\right)\cdot\nabla\psi_i^k \, d\x dt
  = 
  \sum_{t_s=1}^{N_t} \sum_{\x_q=1}^{N_x} \omega_s\omega_q\Delta t\Delta x \F\left(\W(t_s,\x_q)\right)\cdot\nabla\psi_i^k(t_s,\x_q).
\eeq
But the solution is supposed to be regular within the cell (since it is locally approximated by a polynomial), so 
that the Taylor and Lax-Wendroff expansions \eqref{eq:TaylorExpansion} and \eqref{eq:LWExpansion} still apply at any spatial 
quadrature point $\x_q$. 
The same procedure can then be applied at every spatial quadrature point and the whole method is evaluated everywhere and 
each degree of freedom is updated thanks to \eqref{eq:ADERDG}.

For a general algorithm to obtain the successive time derivative as in \eqref{eq:LWExpansion}, 
in particular in the context of the Euler flux, see \cite{DumbserMunz05}
and references therein. 

\subsection{Local space-time predictor ADER procedure}

As we can see in the previous section, the Taylor expansion in time of the solution may lead to extremely complex 
moving from known spatial derivatives to time derivatives through Lax-Wendroff recursive expansion. 
Especially in the context of non-linear fluxes, or even worse, 
when the flux is an unknown black box\dots

Anyway, another much more straightforward method is under development. As far as I know, this method really starts in 
\cite{DumbserEnauxToro08} in the finite volume context and has recent applications in the DG context \cite{Zanotti15,DumbserMunz05} 
(clearly not a throughful bibliography).

The idea here is that, in order to compute the time dependant numerical flux \circled{6} in \eqref{eq:ADERDG}, only 
the outgoing information from each cell is needed. So, on a local space-time mesh, the solution can be evolved within 
each cell independantly. A predictor solution $\Q_k^n(t,x)$ is obtained and is next used to compute integrals \circled{5} and 
\circled{6}. 

Now, we restrict our study to the space-time slab $\T_k^n = \T_k\times[t^n,t^{n+1}]$ and suppose that two functional basis 
$\left(\psi_s^t\right)_{s=1,\dots,N_t}$ and $\left(\psi_q^x\right)_{q=1,\dots,N_x}$ are respectively defined on $[t^n,t^{n+1}]$ and 
$\T_k$, see sections \ref{ssec:DoFs} and \ref{ssec:BasisFunctions}, such that $\Q_k^n$ expands as 
\beq{eq:Qkn}
  \Q_k^n(t,\x) = \sum_{s=1}^{N_t}\sum_{q=1}^{N_x} \Q_k^{s,q} \psi_s^t(t) \psi_q^x(\x).
\eeq
Then $\Q_k^n$ is supposed to locally verify the space-time problem
\beq{eq:XTQ}
  \left\{
    \begin{array}{ll}
      \partial_t \Q_k^n + \Div{\F\left(\Q_k^n\right)} = 0, &\quad \text{ in } \T_k\times]t^n,t^{n+1}],\\
      \Q_k^n(t=0,\x) = \W_k^n(\x), &\quad \text{ in } \T_k,\\
      \Q_k^n(t,\x) = \Q_k^n(t,\x), &\quad \text{ on } \partial\T_k.
    \end{array}
  \right.
\eeq
The last boundary condition is a little bit mysterious and corresponds to a certain personal point-of-view\dots

Next, a local space-time DG formulation is led with an integration by part only in the time direction:
\beq{eq:XTDG}
\begin{split}
  -&\int_{\T_k^n} \Q_k^n(t,\x)\partial_t\psi_s^t(t)\psi_q^x(\x)\,d\x dt\\
  +&\int_{\T_k} \big(\Q_k^n(t^{n+1},\x)\psi_s^t(t^{n+1})\psi_q^x(\x)-\Q_k^n(t^{n},\x)\psi_s^t(t^{n})\psi_q^x(\x)\big)\, d\x\\
  +&\int_{\T_k^n}\Div{\F\left(\Q_k^n\right)}\psi_s^t(t)\psi_q^x(\x)\,d\x dt
  = 0, 
  \quad \forall s\in\saut{1,N_t}, q\in\saut{1,N_x}.
\end{split}
\eeq
If $\Q_k^n$ and $\F(\Q_k^n)$ are supposed to be spanned by the $\psi_s^t\psi_q^x$ space-time basis, and the 
vector of local unknown $\Q_k^{s,q}$ is supposed to be ordered by time layers:
\beq{eq:UnknownOrdering}
  \QQ = \left(\Q_k^{s,q}\right)_{s=1,\dots,N_t,\;q=1,\dots,N_x} 
      = \left(\Q_k^{1,1},\dots,\Q_k^{1,N_x},\dots,\Q_k^{N_t,1},\dots,\Q_k^{N_t,N_x}\right)^T,
\eeq
where superscript $(.)^T$ stands for transposition, then equation \eqref{eq:XTDG} can be rewritten into
\begin{flalign}\label{eq:LocalDGCompact}
  -&\left[\left(\K^t\right)^T\otimes\M^x\right]\QQ 
  + \left[\One_{n+1}^t\otimes\M^x\right]\QQ
  - \left[\One_{n}^t\otimes\M^x\right]\QQ
  + \left[\M^t\otimes\vec{\K}^x\right]\F(\QQ) = 0, \nonumber\\
  \Leftrightarrow&\nonumber\\
  &\left[\left((\K^t)^T-\One_{n+1}^t+\One_n^t\right)\otimes\M^x\right]\QQ = \left[\M^t\otimes\vec{\K}^x\right]\F(\QQ),
\end{flalign}
with the following matrices:
\begin{eqnarray}
  \M_{i,j}^x &=& \int_{\T_k}\psi_i^x(\x)\psi_j^x(\x)\;d\x,\\
  \M_{i,j}^t &=& \int_{t^n}^{t^{n+1}}\psi_i^t(t)\psi_j^t(t)\;dt,\\
  \K_{i,j}^t &=& \int_{t^n}^{t^{n+1}}\psi_i^t(t)\partial_t\psi_j^t(t)\;dt,\\
  \left(\vec{\K}_{i,j}^x\right)_l &=& \int_{\T_k}\psi_i^x(\x)\partial_{x_l}\psi_j^x(\x)\;d\x\\
  (\One_{\xi}^t)_{i,j} &=& \int_{t^n}^{t^{n+1}}\psi_i^t(t^{\xi})\psi_j^t(t)\;dt.
\end{eqnarray}
We have got a local non-linear system which needs to be solved. Fortunately, according to \cite{DumbserBalsaraToro08}, 
the process would be always contractant and the solution can be reached within a few iterations per cells. 

Once this is done, the predictor $\Q_k^n$ is simply used in \circled{5} and \circled{6} 
to update equation \eqref{eq:ADERDG}: 
\begin{align}\label{eq:UpdatePredictor}
 &\quad \forall \T_k \in \Mh, \forall i \in \T_k, \\
  &\M_{ij} \left(\W_j^{k,n+1}-\W_j^{k,n}\right) 
  - 
  \int_{\T_k^n} \F\left(\Q_k^n\right)\cdot\nabla\psi_i^k \, d\x dt 
  + 
  \int_{t^n}^{t^{n+1}}\int_{\partial\T_k}\psi_i^k\FStar(t;\Q^{\text{ext}},\Q^{\text{int}};\n) \, d\s \,dt
  = 
  0. \nonumber
\end{align}

\clearpage
\bibliographystyle{plain}
\bibliography{biblio}

\begin{thebibliography}{10}

\bibitem{Abgrall01}
R.~Abgrall.
\newblock Toward the ultimate conservative scheme: Following the quest.
\newblock {\em J. Comput. Phys}, 167(2):277--315, 2001.

\bibitem{Arnold01}
Douglas~N. Arnold, Franco Brezzi, Bernardo Cockburn, and L.~Donatella Marini.
\newblock Unified analysis of discontinuous galerkin methods for elliptic
  problems.
\newblock {\em SIAM J. Numer. Anal.}, 39(5):1749--1779, May 2001.

\bibitem{BassiRebay2}
F.~Bassi and S.~Rebay.
\newblock A high-order accurate discontinuous finite element method for the
  numerical solution of the compressible navier–stokes equations.
\newblock {\em Journal of Computational Physics}, 131(2):267 -- 279, 1997.

\bibitem{BassiRebay1}
F.~Bassi and S.~Rebay.
\newblock High-order accurate discontinuous finite element solution of the 2d
  euler equations.
\newblock {\em Journal of Computational Physics}, 138(2):251 -- 285, 1997.

\bibitem{CockburnShu2}
B.~Cockburn and Chi-Wang Shu.
\newblock {TVB} {R}unge-{K}utta local projection discontinuous {G}alerkin
  finite element method for conservation laws {II}: General framework.
\newblock {\em Mathematics of Computation}, 52(186):411--435, 1989.

\bibitem{CockburnShu4}
B.~Cockburn and Chi-Wang Shu.
\newblock Tvb runge-kutta local projection discontinuous galerkin finite
  element method for conservation laws iv: the multidimensional case.
\newblock {\em Mathematics of Computation}, 54:545--581, 1990.

\bibitem{CockburnShu1}
B.~Cockburn and Chi-Wang Shu.
\newblock The runge-kutta local projection $p^1$-discontinuous-galerkin finite
  element method for scalar conservation laws.
\newblock {\em M2AN}, 25(3):337--361, 1991.

\bibitem{CockburnShu5}
B.~Cockburn and Chi-Wang Shu.
\newblock The {R}unge-{K}utta discontinuous {G}alerkin method for conservation
  laws {V} - multidimensional systems.
\newblock {\em Journal of Computational Physics}, 141(2):199--124, 1998.

\bibitem{CockburnShu3}
Bernardo Cockburn, San-Yih Lin, and Chi-Wang Shu.
\newblock Tvb runge-kutta local projection discontinuous galerkin finite
  element method for conservation laws iii: One-dimensional systems.
\newblock {\em Journal of Computational Physics}, 84(1):90 -- 113, 1989.

\bibitem{DetteStudden}
H.~Dette and W.~J. Studden.
\newblock {\em The theory of canonical moments with applications in statistics,
  probability, and analysis}.
\newblock John Wiley \& Sons Inc., New York, 1997.

\bibitem{DumbserBalsaraToro08}
Michael Dumbser, Dinshaw~S. Balsara, Eleuterio~F. Toro, and Claus-Dieter Munz.
\newblock A unified framework for the construction of one-step finite volume
  and discontinuous galerkin schemes on unstructured meshes.
\newblock {\em Journal of Computational Physics}, 227(18):8209 -- 8253, 2008.

\bibitem{DumbserEnauxToro08}
Michael Dumbser, Cedric Enaux, and Eleuterio~F. Toro.
\newblock Finite volume schemes of very high order of accuracy for stiff
  hyperbolic balance laws.
\newblock {\em Journal of Computational Physics}, 227(8):3971 -- 4001, 2008.

\bibitem{DumbserMunz05}
Michael Dumbser and Claus-Dieter Munz.
\newblock Building blocks for arbitrary high order discontinuous galerkin
  schemes.
\newblock {\em Journal of Scientific Computing}, 27(1):215--230, 2005.

\bibitem{Gottlieb01}
Sigal Gottlieb, Chi-Wang Shu, and Eitan Tadmor.
\newblock Strong stability-preserving high-order time discretization methods.
\newblock {\em SIAM Review}, 43(1):89--112, 2001.

\bibitem{HesthavenWarburton}
Jan~S. Hestaven and Tim Warburton.
\newblock {\em Nodal Discontinuous Galerkin Methods}.
\newblock Springer-Verlag, 2008.

\bibitem{JohnsonRossmanith13}
Evan~Alexander Johnson and James~A. Rossmanith.
\newblock Outflow positivity limiting for hyperbolic conservation laws. part i:
  Framework and recipe.
\newblock submitted Dec 19, 2012, \url{http://arxiv.org/abs/1212.4695v1}.

\bibitem{OrtnerSuli07}
Christoph Ortner and Endre Süli.
\newblock Discontinuous galerkin finite element approximation of nonlinear
  second-order elliptic and hyperbolic systems.
\newblock {\em SIAM Journal on Numerical Analysis}, 45(4):1370--1397, 2007.

\bibitem{ReedHill73}
W.H. Reed and T.R. Hill.
\newblock Triangular mesh methods for the neutron transport equation.
\newblock In {\em National topical meeting on mathematical models and
  computational techniques for analysis of nuclear systems}, Oct 1973.

\bibitem{ShuSSP88}
Chi-Wang Shu.
\newblock Total-variation-diminishing time discretizations.
\newblock {\em SIAM J. Sci. Stat. Comput.}, 9(6):1073--1084, November 1988.

\bibitem{ShuOsher88}
Chi-Wang Shu and Stanley Osher.
\newblock Efficient implementation of essentially non-oscillatory
  shock-capturing schemes.
\newblock {\em Journal of Computational Physics}, 77(2):439 -- 471, 1988.

\bibitem{ShuOsherII89}
Chi-Wang Shu and Stanley Osher.
\newblock Efficient implementation of essentially non-oscillatory
  shock-capturing schemes, \{II\}.
\newblock {\em Journal of Computational Physics}, 83(1):32 -- 78, 1989.

\bibitem{toro1999}
E.F. Toro.
\newblock {\em Riemann Solvers and Numerical Methods for Fluid Dynamics: A
  Practical Introduction}.
\newblock Springer, New York, 1999.

\bibitem{VanLeer79}
B.~van Leer.
\newblock {Towards the ultimate conservative difference scheme V. A second
  order sequel to Godunov's method}.
\newblock {\em J. Comput. Phys.}, 32:101--136, 1979.

\bibitem{Zanotti15}
Olindo Zanotti, Francesco Fambri, Michael Dumbser, and Arturo Hidalgo.
\newblock Space–time adaptive ader discontinuous galerkin finite element
  schemes with a posteriori sub-cell finite volume limiting.
\newblock {\em Computers \& Fluids}, 118:204 -- 224, 2015.

\bibitem{ZhangThesis}
X.~Zhang.
\newblock {\em Maximum-Principle-Satisfying and Positivity-Preserving High
  Order Schemes for Conservation Laws}.
\newblock PhD thesis, Brown University, 2011.

\bibitem{ZhangShu10a}
Xiangxiong Zhang and Chi-Wang Shu.
\newblock On maximum-principle-satisfying high order schemes for scalar
  conservation laws.
\newblock {\em J. Comput. Phys.}, 229(9):3091--3120, May 2010.

\bibitem{ZhangShu10b}
Xiangxiong Zhang and Chi-Wang Shu.
\newblock On positivity-preserving high order discontinuous galerkin schemes
  for compressible euler equations on rectangular meshes.
\newblock {\em J. Comput. Phys.}, 229(23):8918--8934, November 2010.

\bibitem{ZhangShu11}
Xiangxiong Zhang and Chi-Wang Shu.
\newblock Positivity-preserving high order discontinuous galerkin schemes for
  compressible euler equations with source terms.
\newblock {\em J. Comput. Phys.}, 230(4):1238--1248, February 2011.

\bibitem{ZhangXiaShu12}
Xiangxiong Zhang, Yinhua Xia, and Chi-Wang Shu.
\newblock Maximum-principle-satisfying and positivity-preserving high order
  discontinuous galerkin schemes for conservation laws on triangular meshes.
\newblock {\em J. Sci. Comput.}, 50(1):29--62, January 2012.

\end{thebibliography}

\end{document}